\documentclass{amsart}

\usepackage{amssymb,amsfonts, latexsym,amsmath,hhline,array,longtable}
\usepackage{pdfsync,color, comment,colortbl, mathrsfs,stmaryrd,cite,graphicx,yfonts}

\usepackage{ifpdf}
\ifpdf \usepackage[colorlinks=true, citecolor=blue, linkcolor=blue, urlcolor=blue]{hyperref} \fi

\newtheorem{formula}{}[section]
\newtheorem{definition}[formula]{Definition}
\newtheorem{corollary}[formula]{Corollary}
\newtheorem{remark}[formula]{Remark}
\newtheorem{lemma}[formula]{Lemma}
\newtheorem{theorem}[formula]{Theorem}
\newtheorem*{claim}{Claim}

\def\thrm{\begin{theorem}}
\def\thrml#1{\begin{theorem}\label{#1}}
\def\ethrm{\end{theorem}}
\def\rmrk{\begin{remark}}
\def\rmrkl#1{\begin{remark}\label{#1}}
\def\ermrk{\end{remark}}
\def\dfntn{\begin{definition}}
\def\dfntnl#1{\begin{definition}\label{#1}}
\def\edfntn{\end{definition}}
\def\nmrt{\begin{enumerate}}
\def\enmrt{\end{enumerate}}

\def\qtnl#1{\begin{equation}\label{#1}}
\def\eqtn{\end{equation}}
\def\lmm{\begin{lemma}}
\def\lmml#1{\begin{lemma}\label{#1}}
\def\elmm{\end{lemma}}
\def\crllr{\begin{corollary}}
\def\crllrl#1{\begin{corollary}\label{#1}}
\def\ecrllr{\end{corollary}}
\def\css{\begin{cases}}
\def\ecss{\end{cases}}
\def\prf{\begin{proof}}
\def\eprf{\end{proof}}
\def\clm{\begin{claim}}
\def\eclm{\end{claim}}

\newcommand{\cal}{\mathcal}

\def\cA{{\cal A}}

\def\cS{{\cal S}}

\def\mF{{\mathbb F}}

\def\mZ{{\mathbb Z}}

\DeclareMathOperator{\AGL}{AGL}

\DeclareMathOperator{\cay}{Cay}

\DeclareMathOperator{\Span}{Span}

\DeclareMathOperator{\sym}{Sym}

\makeatletter 
\def\@seccntformat#1{\csname the#1\endcsname. } 
\def\@biblabel#1{#1.} 

\newtheorem{prop}{Proposition}[section]

\newtheorem{corl}[prop]{Corollary}

\theoremstyle{definition}

\def\und#1{{\underline{#1}}}
\def\grp#1{\langle {#1}\rangle}

\begin{document}

\title{Notes on $B$-groups}

\author{Ilia Ponomarenko}
\address{St. Petersburg Department of V.A. Steklov Institute of Mathematics, St. Petersburg, Russia}
\email{inp@pdmi.ras.ru}

\author{Grigory Ryabov}
\address{School of Mathematics and Statistics, Hainan University, Haikou, China}
\address{Novosibirsk State Technical University, Novosibirsk, Russia}
\email{gric2ryabov@gmail.com}

\thanks{The second author was supported by the grant of The Natural Science Foundation of China (project No.~12361003).}
\date{}

\begin{abstract}
Following Wielandt, a finite group $G$ is called a \emph{$B$-group} (\emph{Burnside group}) if every primitive group containing a regular subgroup isomorphic to $G$ is doubly transitive. Using a method of Schur rings, Wielandt proved that every abelian group of composite order which has at least one cyclic Sylow subgroup is a $B$-group. Since then, other infinite families of $B$-groups were found by the same method. A simple analysis of the proofs of these results shows that in all of them a stronger statement was proved for the group $G$ under consideration: every primitive Schur ring over $G$ is trivial. A finite group $G$ possessing the latter property, we call \emph{$BS$-group} (\emph{Burnside-Schur group}). In the present note, we give infinitely many examples of $B$-groups which are not $BS$-groups.
\end{abstract}

\maketitle

\section{$B$-groups}

We assume that the reader is familiar with the basics of permutation groups; definitions of all undefined terms below can be found in  monographs~\cite{Wielandt1964,DixM1996}. The main tools in our proofs are partial difference sets and strongly regular graphs; for relevant definitions and facts, we refer the reader to survey~\cite{MaSL1994}. Throughout the note, cyclic and elementary abelian groups of order~$n$ are denoted by $C_n$ and~$E_n$, respectively.

A group $G$ is defined to be a \emph{$B$-group} if every primitive permutation group containing a regular subgroup isomorphic to $G$ is $2$-transitive. In other words, the $B$-groups are exactly those that  do not appear as regular subgroups of primitive non-$2$-transitive groups; in what follows, the latter groups are called \emph{uniprimitive}.  The concept of a $B$-group was introduced by Wielandt \cite[Definition~25.1]{Wielandt1964} who proved that every abelian group of composite order which has at least one cyclic Sylow subgroup is a $B$-group \cite[Theorem~25.4]{Wielandt1964}. Several other infinite classes of abelian and non-abelian $B$-groups  known at that time were presented   in \cite[Section~25]{Wielandt1964}. It should be noted that the corresponding proofs are mostly based on using the method of Schur rings (see Section~\ref{270224a}).

In~\cite{Cameron1982}, it was observed that the set of natural numbers $n$ for which there exists a primitive group of degree $n$ other than $A_n$ and $S_n$ has zero density in $\mathbb{N}$. This implies that almost all finite groups are $B$-groups, more exactly, all groups of order~$n$ are $B$-groups for almost all positive integers~$n$.  In the sequel, we will need a simple criterion for a group to be a non-$B$-group, see  \cite[Theorem~25.7]{Wielandt1964}. 

\lmml{280224a}
A direct product of at least two groups of the same order at least $3$ is a non-$B$-group.
\elmm

In the 2000s, a new surge of interest in $B$-groups  arose in connection with the completion of the Classification of Finite Simple Groups. The use of the classification is based on the  O'Nan--Scott theorem about the structure of a primitive permutation group~\cite[Theorem~4.1A]{DixM1996}. This technique enables us to study regular subgroups of a primitive group directly, without using Schur rings. In this way, a strong necessary condition was obtained in \cite{LiCH2003} for an abelian group to be a $B$-group. Furthermore, complete classifications were found for the almost simple $B$-groups \cite{Liebeck2010} and for the $B$-groups of order $pq$, where $p$ and $q$ are distinct primes~\cite{Potochnik2001}.

\section{$BS$-groups} \label{270224a}
To explain how Schur rings were used to prove that a group $G$ is a $B$-group, we recall some basic definitions. Let $\cS$ be a partition of $G$ such that $\{1_G\}\in\cS$ and $X^{-1}\in\cS$ for all $X\in\cS$, where $1_G$ is the identity of~$G$ and $X^{-1}=\{x^{-1}:\ x\in X\}$. One can associate with $\cS$ a $\mZ$-submodule
\qtnl{240224a}
\cA=\Span_\mZ\{\und{X}:\ X\in\cS\}
\eqtn
of the group ring $\mZ G$, where the formal sum $\und{X}=\sum_{x\in X}x$ is naturally treated as element of $\mZ G$. The module $\cA$ is called a \emph{Schur ring} over $G$ if $\cA$ is a subring of the ring~$\mZ G$. It was proved by Schur (1933) that  if  a permutation group $K\le \sym(G)$ contains $G$ as a regular subgroup and $\cS$ is a partition of $G$ into the orbits of the  stabilizer of the pont~$1_G$ in~$K$, then the module $\cA=\cA(K,G)$ defined by \eqref{240224a} is a Schur ring over~$G$. 

The Schur ring \eqref{240224a} is said to be \emph{primitive} if for every subgroup $H\le G$, the condition $\und{H}\in\cA$ implies that $H=\{1_G\}$ or $H=G$. An easy example of a primitive Schur ring is given by the \emph{trivial}  Schur ring, i.e., such that $|\cS|=2$. The case $|\cS|=3$ was studied in different contexts. Here we mention only that if an inverse-closed subset $X\subset G$ is not a subgroup of $G$ and $Y$ is the complement to $X$ in the set $G\setminus\{1_G\}$, then $\{1_G,X,Y\}=\cS(\cA)$ for some primitive Schur ring $\cA$ over $G$ if and only if $\grp{X}=G$ and $X$ is a  \emph{regular partial difference set} (\emph{PDS}, for short) in~$G$ in the sense of survey~\cite{MaSL1994}. Note that the latter condition means exactly that the Cayley graph $\cay(G,X)$ over $G$ with connection set $X$ is a connected and coconnected strongly regular graph.

In the above terms, a group $G$ is a $B$-group if the trivial ring is the only primitive Schur ring of the form $\cA(K,G)$, where the  group $K\le \sym(G)$ contains  $G$ as a regular subgroup. A simple analysis shows that in the proofs of all above mentioned results by Wielandt and his colleagues, a more stronger statement than~$G$ is a $B$-group was proved, namely, the only primitive Schur ring over~$G$  is trivial. The groups satisfying the latter property we call \emph{$BS$-groups}. Among them one can find, e.g., dihedral and generalized dicyclic groups. Moreover, it was proved  in~\cite{Bercov1969} that (apart from the minor exception) an abelian group $G$ is not $BS$ if and only if $G$ is the direct product of two subgroups of the same exponent, cf.~Lemma~\ref{280224a}

The above discussion suggests the following natural question, which is the main motivation of the present note.\medskip

{\bf Question.} {\it Is it true that every $B$-group is a $BS$-group?}\medskip

In the next sections, we answer this question negatively by finding several infinite families of $B$-groups  which are not $BS$-groups. 

\section{Exclusive groups}
A finite group $G$ is said to be \emph{exclusive} if $G$ is a $B$-group and $G$ is not a $BS$-group, or equivalently, there is at least one nontrivial primitive Schur ring over $G$ and  the only primitive Schur ring of the form $\cA(K,G)$, where the  group $K\le \sym(G)$ contains  $G$ as a regular subgroup,  is trivial.  The following simple criterion for a group to be a non-$BS$-group immediately follows from the discussion in Section~\ref{270224a}.

\lmml{280224v}
Any group generating by a regular PDS in it is non-$BS$.
\elmm

We make use this criterion to prove that every exclusive group is of order at least~$57$. We will see later, this lower bound on the order of exclusive group is exact.

\lmml{upto56}
Each group of order at most~$56$ is non-exclusive.
\elmm

\begin{proof}
The list of all non-$BS$-groups  of order at most~$63$ was computed by Ziv-Av, see~\cite[Table~3]{Ziv-av2014}.\footnote{It should be noted that  the term $B$-group in that paper means the same as what we call $BS$-group in the present paper.}  It suffices to verify that each group of this list that has order $\le 56$ is not a $B$-group. Exactly $10$ groups of Ziv-Av's list are direct products of at least two groups of the same order at least~$3$. Thus each of these groups is a  non-$B$-group (Lemma~\ref{280224a}). The remaining groups are given in Table~\ref{tbl}.
\begin{table}[t]
\begin{tabular}{|c|c|c|}
\hline
No. & Order        &  Groups    \\
\hline
$1.$    & $16$     &   $(C_4\times C_2):C_2$, $C_4:C_4$, $C_8:C_2$, $QD_{16}$,  $ C_2\times D_8$\\
\hline
$2.$	& $21$    &  $C_7:C_3$   \\
\hline
$3.$	& $27$    &  $E_9:C_3$, $C_9:C_3$   \\
\hline
$4.$	& $36$    &  $C_3\times (C_3:C_4)$, $E_9:C_4$, $C_3\times A_4$, $C_2\times (E_9: C_2)$ \\
\hline
$5.$	& $55$    &  $C_{11}:C_5$   \\
\hline
$6.$	& $56$    &  $\AGL_1(8)$   \\
\hline
\end{tabular}
\caption{Non $BS$-groups of order $\le 56$.}
\label{tbl}
\end{table}

The fact that the two nonabelian groups of order~$27$ in the third row of Table~\ref{tbl} are non-$B$-groups follows from, e.g.,~\cite[p.~55]{Jones1979}. The two nonabelian groups of orders~$21$ and $55$ in the second and fifth rows of Table~\ref{tbl} are not $B$-groups by ~\cite[Theorem~1.1]{Potochnik2001}. The group $\AGL_1(8)$ of order~$56$ in the six row of Table~\ref{tbl} is  a non-$B$-group, because it appears in~\cite[Table~16.2]{Liebeck2010} as a regular subgroup of the alternating group~$A_8$ in its uniprimitive action on triples. 
	
It still remain  to consider the groups of orders $16$ or $36$ in the first and fourth rows of Table~\ref{tbl}. Each of them appears as a regular automorphism group of the  lattice graphs $L_2(4)$, $L_2(6)$, $L_3(C_6)$, see~\cite[Tables~4.2,~4.4,~4.5]{Heinze2001}. Since the automorphism group of each of this graph is uniprimitive, we conclude that each group in the first or fourth row of Table~\ref{tbl} is a non-$B$-group.
\end{proof}

Now we construct an infinite families of exclusive groups, which contains a minimal example of exclusive group.

\thrml{190224a}
Let $p\equiv 7\mod~12$ be a prime, $p\ne 7, 31$. Then the Frobenius group $C_p:C_{(p-1)/6}$ is exclusive.
\ethrm

\prf
Let $\mF$ be a finite field of order $p$, and let $A$ be the  subgroup of index $6$ in the multiplicative group of~$\mF$. We are interested in a subgroup 
$$
G=\{x\mapsto ax+b,~x\in \mF:\ a\in A,\ b\in \mF\}.
$$
of the 
group $\AGL_1(p)$. Clearly, $|G|=p(p-1)/6$. In accordance with \cite[Theorem~3.2(2b)]{Clapham1976}, there exists  an element $x\in \mF$  such that the set
$$
B=\{\{{0_\mF}^g,{1_\mF}^g,x^g\}:~g\in G\},
$$
forms a Steiner triple system on the elements of  $G$; in particular, $B$ consists of exactly $p(p-1)/6=|G|$ triples and the permutation group $G\le\sym(B)$ is regular.  

Denote by $\Gamma$ the \emph{block} graph of the Steiner triple system $B$: it has vertex set~$B$ and two vertices  are adjacent  if and only if the corresponding blocks intersect nontrivially. It is well known that the block graph of any Steiner triple system, in particular $\Gamma$, is strongly regular, connected, and coconnected. Since $G$ is a regular automorphism group of $\Gamma$, we conclude that $\Gamma\cong\cay(G,X)$ for some $X\subset G$ such that $\grp{X}=G$. It follows (see Section~\ref{270224a}) that $X$ is a regular PDS that generates $G$. Thus $G$ is a non-$BS$-group by Lemma~\ref{280224v}.

It remains to prove that $G$ is a $B$-group. To this end, we need to summarize some simple facts about the order of~$G$.

\clm
The number $n=|G|$ is odd and composite. Moreover, $n=a^b$ for no $a\ge 2$ and $b\ge 2$.
\eclm
\prf
The first statement follows from the equality $n=p(p-1)/6$ and the choice of~$p$. Since $n$ is divisible by $p$ and not by $p^2$, then in any representation  $n=a^b$ with $a\ge 2$, the integer $b$ equals~$1$. This proves the second statement.
\eprf

Now let $K$ be a uniprimitive group containing $G$ as a regular subgroup. Then the degree $n$ of this group is as  in the Claim. By the O'Nan--Scott theorem~\cite[Theorem~4.1A]{DixM1996}, this implies that $K$ is an almost simple group. At this point, we will use a complete classification of all pairs $(K',G')$, where $K'$ is a primitive almost simple group and $G'$ is a regular subgroup of~$K'$, see paper~\cite{Liebeck2010}. Inspecting [Tables~16.1--16.3] of that paper and discarding $2$-transitive groups $K'$ and groups~$G'$ of even order, we arrive at the two possibilities: either
\qtnl{020324a}
n=q\cdot\frac{(q-1)}{2},\ q \text{ is a prime power,}\ q\equiv 3 \pmod 4,
\eqtn
or 
$$
n\in \{31\cdot 5,\ 11\cdot 5,\ 23\cdot 11,\  29\cdot 7,\ 59\cdot 29,\ 3\cdot 19\cdot 9,\ 15\}.
$$ 
In the latter case, $n=31\cdot 5$ is not possible, because $p\ne 31$, and the other values of~$n$ are exculded by the condition $p\equiv 7\mod~12$. Finally, assume that relations~\eqref{020324a} hold. Then $p(p-1)/6=n=q(q-1)/2$. It follows that $p$ divides $q$, for  otherwise, $p$ divides $(q-1)/2$ and 
$$
(p-1)/6 <p\le (q-1)/2<q,
$$
whence $p(p-1)/6<q(q-1)/2=n$, a contradiction. A similar argument shows that $p$ cannot divide $q$. This implies that $G$ cannot be a regular subgroup of $K$, a contradiction. Thus $G$ is a $B$-group, which completes the proof.
\eprf

The smallest prime $p$ satisfying the condition Theorem~\ref{190224a} is equal to~$19$. Together with Lemma~\ref{upto56}, this  proves the following statement.

\begin{corl}\label{57}
$C_{19}: C_3$ is an exclusive group of the smallest possible order.
\end{corl}

\section{More examples of exclusive groups}\label{96}

Inspecting the list of all small primitive groups in~\cite{Roney2005}, one can see that there are no  uniprimitive permutation groups of degree~$96$. It follows that any group of order $96$ is a $B$-group. On the other hand, according to the computations made in~\cite{Golemac2006}, there are  nine nonabelian groups of order~$96$ that admit a regular PDS. Thus all these groups are exclusive by Lemma~\ref{280224v}.

Let $G$  be an extraspecial group of order~$p^3$ and exponent~$p^2$. It was observed in \cite[p.~55]{Jones1979} that $G$ is a $B$-group if $p\geq 5$. On the other hand, the group $G$ contains a regular PDS of cardinality  $p^2+p-2$ for all $p\ge 3$, see ~\cite[Theorem~1]{Swartz2015}. Since this PDS generates $G$, we conclude by Lemma~\ref{280224v} that $G$ is an exclusive group for all~$p\ge 5$. 

\thrm
An extraspecial group of order~$p^3$ and exponent~$p^2$ is an exclusive group for all primes $p\geq 5$.
\ethrm

In almost the same technique as used in the proof of Theorem~\ref{190224a}, one can find one more infinite family of exclusive groups. Namely, from~\cite[Theorem~12.16]{MaSL1994}, it follows that a group 
$$
G=E_q\times (E_r:C_{(r-1)/2q})
$$ 
admits a regular PDS if $q$ is a $2$-power, $r\ge 9$ is a $3$-power, and $2q$ strictly divides $r-1$. Under additional assumptions on $q$ and $r$, one can verify (using again a classification of regular subgroups in an almost simple primitive group~\cite[Tables~16.1-16.3]{Liebeck2010}) that $G$ is a $B$-group for infinitely many~$q$ and~$r$.

\hspace{5mm}

The authors would like to thank Prof. M. Grechkoseeva and Prof. A. Vasil'ev for the discussions on the subject matters. 


\providecommand{\bysame}{\leavevmode\hbox to3em{\hrulefill}\thinspace}
\providecommand{\MR}{\relax\ifhmode\unskip\space\fi MR }
\providecommand{\MRhref}[2]{%
	\href{http://www.ams.org/mathscinet-getitem?mr=#1}{#2}
}
\providecommand{\href}[2]{#2}

\end{document}